\documentclass{article}

\usepackage{amsmath,amsthm}

\newtheorem{theorem}{Theorem}

\newtheorem{definition}[theorem]{Definition}

\newtheorem{comment}[theorem]{Comment}

\newtheorem{lemma}[theorem]{Lemma}

\newcommand{\dR}{{{\bf R}}}
\newcommand{\ep}{\varepsilon}
\newcommand{\dN}{{{\bf N}}}

\begin{document}

\title{Borel Games with Lower-Semi-Continuous Payoffs%
\thanks{This research was supported by the Israel Science Foundation (grant number 212/09).}}

\author{Ayala Mashiah-Yaakovi\thanks{School of Mathematical Sciences, Tel Aviv University, Tel
Aviv 69978, Israel. e-mail: \textsf{ayalam@post.tau.ac.il}.}
and Eilon Solan\thanks{School of Mathematical Sciences,
Tel Aviv University, Tel Aviv 69978, Israel. e-mail:
\textsf{eilons@post.tau.ac.il}}}

\maketitle

\begin{abstract}
We prove that every two-player non-zero-sum Borel game with
lower-semi-continuous payoffs admits a subgame-perfect $\ep$-equilibrium.
This result complements Example 3 in Solan and Vieille (2003),
which shows that a subgame-perfect $\ep$-equilibrium need not exists when the payoffs are not lower-semi-continuous.
\end{abstract}

\section{Introduction}

Borel games are sequential games where two players alternately choose actions.
The payoff of each player is a function of the infinite sequence of actions that the players chose.
Borel games were introduced by Gale and Stuart (1953), who studied zero-sum games
where the payoff function is the indicator of some set.
In other words, player 1 wins if the play generated by the players is in a given set of plays,
and player 2 wins otherwise.
Martin (1975) proved that if the winning set of player 1 is Borel measurable, then the game is determined:
either player 1 has a winning strategy or player 2 has a winning strategy.
This result implies that every two-player zero-sum Borel game has a value,
provided the payoff function is bounded and measurable.

Mertens and Neyman (see Mertens, 1987) used the existence of the value in multi-player non-zero-sum Borel games
to prove that for every $\ep > 0$, every
multi-player non-zero-sum Borel game has an $\ep$-equilibrium,
provided the payoff functions are bounded and measurable.
The $\ep$-equilibrium strategies constructed by Mertens and Neyman are as follows:
each player $i$ starts by following an $\frac{\ep}{2}$-optimal strategy in an auxiliary two-player zero-sum game $G_i$,
where the payoff is that of player $i$, player $i$ is the maximizer and the other players try to minimize player $i$'s
payoff.
This goes on as long as no player deviates.
Once some player, say player $i$, deviates,
the other players switch to an $\frac{\ep}{2}$-optimal strategy of the minimizers in the game $G_i$.

Thus, the players start by generating a play that yields all of them a high payoff,
and, if a player deviates, he is punished with a low payoff.
This construction has the disadvantage that in the punishment phase,
the punishers may lower their own payoffs.
Therefore, in real-life situations, players may be reluctant to follow the equilibrium strategies
constructed by Mertens and Neyman.

To deal with such non-credible threats of punishment,
Selten (1965, 1973) introduced the concept of subgame-perfect equilibrium.
A strategy vector is a subgame-perfect $\ep$-equilibrium if it induces an $\ep$-equilibrium
after any possible finite history of actions.
Ummels (2005) proved the existence of a subgame-perfect equilibrium when
the payoff function of each player is the indicator of some set.
His proof is based on the following recursive construction.
First, one identifies all finite histories which are a winning position to either player 1 or player 2;
that is, if this finite history occurs, one of the players can ensure that his payoff is 1.
After each such finite history one instruct the winning player to play his winning strategy,
and the other player is instructed to play his best response.
One then identifies winning positions to the two players in the new game,
assuming the behavior in the first set of winning positions is set.
The process repeats itself, until it reaches a stable state.
Ummels (2005) proves that the pair of strategies thus defined is a subgame-perfect equilibrium.

In the present paper we show that every Borel game with bounded and lower-semi-continuous payoffs admits
a subgame-perfect $\ep$-equilibrium, for every $\ep > 0$.
This result complements Example 3 in Solan and Vieille (2003) that shows that in games with general payoff
functions a subgame-perfect $\ep$-equilibrium
need not exist.

The determinacy of Borel games has attracted a lot of attention in
descriptive set theory (see, e.g., Schilling and Vaught (1983) and Kechris (1995)).
A rich literature identifies winning positions for the two players in
the class of games that are played on graphs
(see Gr\"adel (2004) for a survey).
Two-player zero-sum Borel games where used in the computer science literature
to study reactive non-terminating programs (see, e.g., Thomas (2002))
and model checking in $\mu$-calculus (see, e.g., Emerson et al. (2001)),
and in economics to show that measurable tests are manipulable (Shmaya, 2008).

Our result also relates to the game theoretic literature that
studies the existence of a subgame-perfect $\ep$-equilibrium in various classes of infinite games,
see, e.g., Solan and Vieille (2003), Mashiah-Yaakovi (2009) or Flesch et al. (2008, 2009).
In particular, our result generalizes results in Flesch et al. (2008, 2009).

The paper is organized as follows.
The model and the main result appear in Section \ref{model}.
Section \ref{proof} contains the proof of the main result,
and Section \ref{comments} contains some comments.

\section{The Model and the Main Result}
\label{model}

\begin{definition}
A (two-player non-zero-sum) {\em Borel game} is a triplet $(A, u^1,u^2)$ where
$A$ is a set of actions, and
$u^1,u^2 : A^\dN \to \dR$ are payoff functions.
\end{definition}

The game is played as follows.
At every odd stage $t$ player 1 chooses an action $a_t \in A$,
and at every even stage $t$ player 2 chooses an action $a_t \in A$.
While making his choice at stage $t$, the player knows the sequence $(a_j)_{j < t}$
of actions that was chosen by the players in previous stages.
The payoff to each player $i$ is $u^i(a_1,a_2,\ldots)$.
The description of the game is common knowledge among the players.
\begin{comment}
The assumptions that (1) there are only two players,
(2) both players have the same action set,
and (3) the set of action is the same in all stages,
are made for simplicity of notations only.
Nothing that is said below would be affected if there are more than two players,
if the players have different action sets,
or if the sets of actions depend on past choices of actions.
\end{comment}

Denote by $\emptyset$ the empty history at the beginning of the game.
The set of possible histories%
\footnote{By convention, $A^0 = \{\emptyset\}$ contains only the empty history.}
at stage $t$ is $H_t = A^{t-1}$.
Denote by $H^1 = \bigcup_{k \in \dN} H_{2k-1}$ and
$H^2 = \bigcup_{k \in \dN} H_{2k}$ the sets of possible histories at decision points
of player 1 and player 2 respectively.
\begin{definition}
A {\em strategy} for player $i$ is a function $\sigma^i : H^i \to A$.
\end{definition}
We denote by $\Sigma^i$ the strategy space of player $i$.
Every pair of strategies $(\sigma^1,\sigma^2) \in \Sigma^1 \times \Sigma^2$ determines a unique play
$p = p(\sigma^1,\sigma^2) = (a_t)_{t \in \dN} \in A^\dN$ as follows:
\begin{eqnarray}
a_1 &=& \sigma^1(\emptyset),\\
a_{2k} &=& \sigma^2(a_1,a_2,\ldots,a_{2k-1}), \ \ \ \forall k \in \dN,\\
a_{2k+1} &=& \sigma^1(a_1,a_2,\ldots,a_{2k}), \ \ \ \forall k \in \dN.
\end{eqnarray}
We denote by $u^i(\sigma^1,\sigma^2) = u^i(p(\sigma^1,\sigma^2))$
the payoff of player $i$ when the two players follow $(\sigma^1,\sigma^2)$.

\begin{definition}
Let $\ep \geq 0$.
A pair of strategies $(\sigma^{1}_*,\sigma^2_*)$ is an {\em $\ep$-equilibrium}
if
\begin{eqnarray}
u^1(\sigma_*^1,\sigma_*^2) &\geq& u^1(\sigma^1,\sigma_*^2)-\ep, \ \ \ \forall \sigma^1 \in \Sigma^1,\\
u^2(\sigma_*^1,\sigma_*^2) &\geq& u^2(\sigma_*^1,\sigma^2)-\ep, \ \ \ \forall \sigma^2 \in \Sigma^2.
\end{eqnarray}
\end{definition}

Throughout the paper we endow $A$ with the discrete topology,
and $A^\dN$ with the product topology.

The game is called {\em zero-sum} if $u^1(p) + u^2(p) = 0$ for every $p \in A^\dN$.
The result of Martin (1975) implies that
in zero-sum games, an $\ep$-equilibrium exists for every $\ep > 0$ under merely a measurability condition.
\begin{theorem}
If the game is zero-sum, and if $u^1$ is bounded and Borel measurable,
then an $\ep$-equilibrium exists for every $\ep > 0$.
\end{theorem}
This results implies the existence of an $\ep$-equilibrium in every two-player non-zero-sum game.
\begin{theorem}[Mertens, 1987]
If $u^1$ and $u^2$ are bounded and Borel measurable,
then an $\ep$-equilibrium exists for every $\ep > 0$.
\end{theorem}

A stronger notion of equilibrium is the notion of subgame-perfect equilibrium.
Every finite history $h = (a_1,a_2,\ldots,a_l) \in H^1 \cup H^2$,
together with a pair of strategies $(\sigma^1,\sigma^2)$,
determines an infinite play $p = p(\sigma^1,\sigma^2 \mid h) = (b_t)_{t \in \dN} \in A^\dN$ as follows:
\begin{eqnarray}
b_j &=& a_j, \ \ \ \ \ 1 \leq j \leq l,\\
b_{2j} &=& \sigma^{2}(b_1,b_2,\ldots,b_{2j-1}), \ \ \ l < 2j,\\
b_{2j+1} &=& \sigma^{1}(b_1,b_2,\ldots,b_{2j}), \ \ \ l < 2j+1.
\end{eqnarray}
We denote by $u^i(\sigma^1,\sigma^2 \mid h) = u^i(p(\sigma^1,\sigma^2 \mid h))$ the payoff of player $i$ at this play.

\begin{definition}
Let $\ep \geq 0$.
A pair of strategies $(\sigma^{1}_*,\sigma^2_*)$ is a {\em subgame-perfect $\ep$-equilibrium}
if for every finite history $h \in H^1 \cup H^2$ one has:
\begin{eqnarray}
u^1(\sigma_*^1,\sigma_*^2 \mid h) &\geq& u^1(\sigma^1,\sigma_*^2 \mid h)-\ep, \ \ \ \forall \sigma^1 \in \Sigma^1,\\
u^2(\sigma_*^1,\sigma_*^2 \mid h) &\geq& u^2(\sigma_*^1,\sigma^2 \mid h)-\ep, \ \ \ \forall \sigma^2 \in \Sigma^2.
\end{eqnarray}
\end{definition}

Thus, every finite history $h$ defines a subgame ---
the subgame that starts once the finite history $h$ occurs.
A strategy pair is a subgame-perfect $\ep$-equilibrium if it induces an $\ep$-equilibrium in all subgames.

We say that a finite history $h = (a_t)_{t=1}^l$ is a {\em prefix} of the play $p = (b_t)_{t \in \dN} \in A^\dN$,
or that $p$ is an {\em extension} of $h$,
if $a_t=b_t$ for every $t \in \{1,2,\ldots,l\}$, and we denote it by $h \prec p$.
We say that a finite history $h = (a_t)_{t=1}^l$ is a {\em prefix} of the finite history $h' = (b_t)_{t=1}^m \in A^\dN$,
or that $h'$ is an {\em extension} of $h$,
if $l \leq m$ and $a_t=b_t$ for every $t \in \{1,2,\ldots,l\}$, and we denote it by $h \preceq h'$.

When $A$ is endowed with the discrete topology, and $A^\dN$ is endowed with the product topology,
then a sequence $(p^k)_{k \in \dN}$ of plays converges to a limit $p$
if and only if every prefix $h$ of $p$ is a prefix of all the plays $(p^k)_{k \in \dN}$
except possibly of finitely many of them.

\begin{definition}
The payoff function $u^i$ is {\em lower-semi-continuous}
if for every sequence $(p^k)_{k \in \dN}$ of infinite plays in $H^1 \cup H^2$ that converges to a limit $p$
one has
\[ \liminf_{k \to \infty} u^i(p^k) \geq u^i(p). \]
\end{definition}

Note that every lower-semi-continuous function is Borel measurable.
Our main result is the following.

\begin{theorem}
\label{theorem5}
If $u^1$ and $u^2$ are lower-semi-continuous and bounded,
then the game admits a subgame-perfect $\ep$-equilibrium for every $\ep > 0$.
\end{theorem}

\section{Proof}
\label{proof}

We first note that if one changes the payoffs $u^1$ and $u^2$ by at most $\ep$
(in the supremum norm),
then every subgame-perfect $\ep$-equilibrium in the original game is a subgame-perfect $3\ep$-equilibrium
in the new game.
Because the range of the payoff functions is bounded,
we can assume w.l.o.g. that the range of $u^1$ and $u^2$ is in fact finite.
We will show that if the payoff functions are bounded and have finite range, there is
a subgame-perfect 0-equilibrium.

For $i \in \{1,2\}$ we denote by $-i$ the player who is not $i$.
We denote by $i(h)$ the player who has to choose an action after the history $h$;
$i(h) = 1$ if $h$ has even length, and $i(h) = 2$ otherwise.
We denote by $-i(h)$ the player who is not $i(h)$.

\subsection{Subgame-perfect optimal strategies}

\begin{definition}
Let $h \in H^1 \cup H^2$ be a finite history.
The real number $v(h)$ is called the {\em value at $h$} if
\begin{equation}
\label{equ 101}
v(h) =
\max_{\sigma^{i(h)} \in \Sigma^{i(h)}} \min_{\sigma^{-i(h)} \in \Sigma^{-i(h)}} u^{i(h)}(\sigma^1,\sigma^2 \mid h)
= \min_{\sigma^{-i(h)} \in \Sigma^{-i(h)}} \max_{\sigma^{i(h)} \in \Sigma^{i(h)}} u^{i(h)}(\sigma^1,\sigma^2 \mid h).
\end{equation}
\end{definition}
Because the range of the functions $u^1$ and $u^2$ is finite,
the maximum and minimum in (\ref{equ 101}) are well defined.

\begin{definition}
A strategy $\sigma^i_* \in \Sigma^i$ of player $i$ is called
{\em subgame-perfect optimal} if for every $h \in H^i$  one has
\[ u^i(\sigma^i_*,\sigma^{-i} \mid h) \geq v(h), \ \ \ \forall \sigma^{-i} \in \Sigma^{-i}. \]
\end{definition}

The next result follows from Martin (1975).
\begin{theorem}
\label{theorem01}
In every zero-sum game with bounded and discrete payoffs,
if $u^1$ is Borel measurable then
both players have subgame-perfect optimal strategies.
\end{theorem}

\subsection{Constructing some sequences}

In this subsection we define for every finite history $h \in H^1 \cup H^2$
and every ordinal $\xi$, (a) a real number $\alpha_\xi(h)$,
and (b) a set $H_\xi(h)$ of plays.
The sequence $(\alpha_\xi(h))_\xi$ will be a non-decreasing sequence of lower bounds
to subgame-perfect 0-equilibrium payoff for the player who makes the decision at $h$,
in the game that start at $h$.
The sequence $(H_\xi(h))_\xi$ will be a non-increasing (by inclusion) sequence of sets of histories,
such that all plays generated by a subgame-perfect 0-equilibrium at the game that starts at $h$
are included in every set in this sequence.

Every play $p = (a_t)_{t \in \dN} \in A^\dN$ defines a sequence $(h_n)_{n \in \dN}$ of finite histories,
where $h_n = (a_t)_{t=1}^n$ is the prefix of length $n$ of $p$.

Suppose that one is given a real-valued function $\alpha$ that is defined over the set $H^1 \cup H^2$
of finite histories.

\begin{definition}
Let $h$ be a finite history, and $p$ a play that extends $h$.
The play $p$ is called {\em $\alpha$-monotonic} at $h$ if
the two sequences $(\alpha(h_{2k}))_{\{k \colon h \preceq h_{2k}\}}$ and
$(\alpha(h_{2k+1}))_{\{k \colon h \preceq h_{2k+1}\}}$
are non-decreasing.
\end{definition}

A play $p$ is called {\em $\alpha$-viable} at a given finite history $h$
if for each player $i$,
the payoff $u^i(p)$ that $p$ yields to player $i$ is higher than $\alpha(h')$,
for every prefix $h'$ of $p$ that extends $h$ and that is a decision history for player $i$.
Formally,\footnote{This definition is adapted from Flesch et al. (2008).}
\begin{definition}
Let $\alpha : H^1 \cup H^2 \to \dR$ be a real-valued function,
and let $h$ be a finite history.
We say that a play $p$ is {\em $\alpha$-viable} at $h$ if
\begin{itemize}
\item   $h$ is a prefix of $p$,
\item   $u^{i(h')}(p) \geq \alpha(h')$ for every history $h' \in H^1 \cup H^2$ that satisfies $h \preceq h' \prec p$.
\end{itemize}
\end{definition}

The following lemma lists two simple properties of the play that is realized when both players follow
subgame-perfect optimal strategies.
It follows from the definitions, and its proof is omitted.
\begin{lemma}
\label{lemma simple1}
Fix a finite history $h \in H^1 \cup H^2$.
Let $\sigma^1$ and $\sigma^2$ be subgame-perfect optimal strategies of the two players in the subgame
that starts at $h$.
Then
the realized play $p(\sigma^1,\sigma^2)$ is both $\alpha_1$-monotonic and $\alpha_1$-viable at $h$.
\end{lemma}

Set
\begin{eqnarray}
\label{equ101}
\alpha_{1}(h) &=&
v^{i(h)}(h),\\
H_1(h) &=&
\{ p \in A^\dN \colon h \prec p, p \hbox{\ is\ } \alpha_1\hbox{-viable at } h\}.
\label{equ101a}
\end{eqnarray}
Thus, $\alpha_1(h)$ is the value of the zero-sum subgame that ``starts at $h$''
with payoffs that are the payoffs of player $i(h)$.
This is a lower bound on the equilibrium payoff of player $i(h)$ in this subgame.
$H_1(h)$ is the set of all plays that yield both players at least their value in all subgames that may occur after $h$.

If $h = (a_t)_{t=1}^l$ is a finite history with length $l$, and $a \in A$,
we denote by $(h,a) = (a_1,a_2,\ldots,a_l,a)$ the history of length $l+1$ that starts with $h$ and ends with $a$.

For every successor ordinal $\xi+1$, define
\begin{eqnarray}
\label{equ102}
\alpha_{\xi+1}(h) &=&
\max_{a \in A} \min_{p \in H_\xi(h,a)}u^{i(h)}(p),\\
H_{\xi+1}(h) &=&
\left\{ p \in \cup_{a \in A} H_\xi(h,a) \colon u^{i(h)}(p) \geq \alpha_{\xi+1}(h)\right\}.
\label{equ102a}
\end{eqnarray}
If $H_\xi(h,a)$ contains all plays that can be generated by subgame-perfect 0-equilibria
in the subgame that starts at $(h,a)$,
then in every subgame-perfect 0-equilibrium of the subgame that starts at $h$
player $i(h)$ will receive at least $\alpha_{\xi+1}(h)$.
Moreover, after the history $h$, player $i(h)$ will not play an action that does not maximize
the right-hand side of (\ref{equ102}),
and therefore $H_{\xi+1}(h)$ contains all plays that can be generated by subgame-perfect 0-equilibria
in the subgame that starts at $h$.

For a limit ordinal $\xi$ define
\begin{eqnarray}
\label{equ103}
\alpha_\xi(h) &=& \max_{\lambda < \xi} \alpha_\lambda(h),\\
H_\xi(h) &=& \{ p \in A^\dN \colon h \preceq p, p \hbox{ is } \alpha_\xi\hbox{-viable at } h\}.
\label{equ103a}
\end{eqnarray}

The following observation, which follows from the definitions, will be used later.
\begin{lemma}
\label{lemma simple0}
Suppose that $p \in H_\xi(h,a)$,
where $a \in A$ achieves the maximum in the right-hand side of (\ref{equ102}).
Then $p \in H_{\xi+1}(h)$.
\end{lemma}

\subsection{Properties of the sequences $(\alpha_\xi(h))_\xi$ and $(H_\xi(h))_\xi$}

The following theorem states some properties of the sequences $(\alpha_\xi(h))_\xi$ and $(H_\xi(h))_\xi$,
which play a crucial role in the proof of the main result.
\begin{theorem}
\label{theorem17}
The following holds for every $h \in H^1 \cup H^2$:
\begin{enumerate}
\item
The set $H_\xi(h)$ is not empty for every ordinal $\xi$.
\item
The sequence $(H_\xi(h))_\xi$ is monotonic non-increasing (by inclusion).
\item
The sequence $(\alpha_\xi(h))_\xi$ is monotonic non-decreasing.
\item
For $\xi=1$ and for every limit ordinal $\xi$,
there is a play $p \in H_\xi(h)$ that is $\alpha_\xi$-monotonic at $h$.
\end{enumerate}
\end{theorem}

\begin{proof}
The proof is by transfinite induction.

\noindent\underline{Part 1:}
$H_1(h) \neq \emptyset$ for every $h \in H^1 \cup H^2$.
Moreover, there is $p \in H_1(h)$ that is $\alpha_1$-monotonic at $h$.

This fact follows from Theorem \ref{theorem01} and Lemma \ref{lemma simple1}:
once both players follows a subgame-perfect optimal strategy,
the realized play is $\alpha_1$-viable and $\alpha_1$-monotonic at $h$.

\noindent\underline{Part 2:}
$\alpha_2(h) \geq \alpha_1(h)$ for every $h \in H^1 \cup H^2$.

Every play in $H_1(h,a)$ is $\alpha_1$-viable at $(h,a)$.
Therefore $\alpha_2(h)$ is the value of the game where at the second stage player $i(h)$
receives the minimum payoff generated by an $\alpha_1$-viable plays.
But $\alpha_1(h)$ is the value of the game where at the second stage player $i(h)$
receives the minimum payoff generated by a play which satisfies the $\alpha_1$-viability condition
only for histories that are controlled by player $i(h)$.
This implies that indeed $\alpha_2(h) \geq \alpha_1(h)$.

\noindent\underline{Part 3:}
$H_2(h) \subseteq H_1(h)$ for every $h \in H^1 \cup H^2$.

This follows from Part 2 and Eq. (\ref{equ102a}).

\noindent\underline{Part 4:}
If
$\alpha_{\xi+1}(h) \geq \alpha_\xi(h)$ and $H_{\xi+1}(h) \subseteq H_\xi(h)$
for every $h \in H^1 \cup H^2$,
then
$\alpha_{\xi+2}(h) \geq \alpha_{\xi+1}(h)$ and $H_{\xi+2}(h) \subseteq H_{\xi+1}(h)$
for every $h \in H^1 \cup H^2$.

This follows from the definitions (\ref{equ102}) and (\ref{equ102a}).

\noindent\underline{Part 5:}
For every limit ordinal $\xi$, every ordinal $\lambda < \xi$, and every $h \in H^1 \cup H^2$, one has
$\alpha_{\xi}(h) \geq \alpha_\lambda(h)$ and $H_{\xi}(h) \subseteq H_\lambda(h)$.

This follows from the definitions (\ref{equ103}) and (\ref{equ103a}).

\noindent\underline{Part 6:}
For every ordinal $\xi$ and every $h \in H^1 \cup H^2$,
If $H_{\xi}(h) \neq \emptyset$ then
$H_{\xi+1}(h) \neq \emptyset$.

This follows from the definitions (\ref{equ102}) and (\ref{equ102a}).

\noindent\underline{Part 7:}
For every limit ordinal $\xi$ and every $h \in H^1 \cup H^2$ one has
$H_{\xi}(h) \neq \emptyset$.
Moreover, there is a play $p \in H_\xi(h)$ that is $\alpha_\xi$-monotonic at $h$.

This is the difficult part of the proof.
Fix a limit ordinal $\xi$ and a finite history $h$.
We are going to generate a play that extends $h$, and we will show that it is in $H_\xi(h)$
and it is $\alpha_\xi$-monotonic at $h$.
The play will be generated in iterations,
where the construction in odd iterations differs from the construction in even iterations.
We will then prove that an infinite play is generated after an even number of iterations.
Finally we will prove that it is in $H_\xi(h)$ and that it is $\alpha_\xi$-monotonic at $h$.

\noindent\textbf{Odd iterations:}

Let $h_1$ be the history at the beginning of the iterations.
For the first iteration, $h_1=h$.
For all other odd iterations, it is the history generated by the previous even iteration.

Consider the following algorithm that generates a finite history or a play that extends $h_1$.
\begin{enumerate}
\item   Let $\xi_1 < \xi$ be a successor ordinal that satisfies $\alpha_\xi(h_1) = \alpha_{\xi_1}(h_1)$.
Such an ordinal exists because every set of ordinals has a minimal element.
\item   Let $a_1$ be an action of player $i(h_1)$ that achieves the maximum in (\ref{equ102}) for $h_1$ and $\xi_1$.
\item   Set $h_2 = (h_1,a_1)$.
\item   Let $\xi_2 \geq \xi_1-1$ be the minimal ordinal that satisfies $\alpha_\xi(h_2) = \alpha_{\xi_2}(h_2)$.
Note that because $\xi$ is a limit ordinal, $\xi_2 < \xi$.
%Note that $\xi_2 < \xi$, and $\xi_2$ is a successor ordinal provided $\xi_2 > 1$.
%Indeed, because the range of the payoffs is finite,
%for every limit ordinal $\lambda > 1$
%there is an ordinal $\mu < \lambda$ such that $\alpha_mu(h_2) = \alpha_\lambda(h_2)$,
%and therefore $\xi_2$ cannot be a limit ordinal larger than 1.
\item   Let $a_2$ be an action of player $i(h_2)$ that achieves the maximum in (\ref{equ102}) for $h_2$ and $\xi_2$.
\item   Set $h_3 = (h_1,a_1,a_2)$.
\item   Continue this way to create a sequence $(h_1,\xi_1,a_1,h_2,\xi_2,a_2,\ldots)$ that extends $h_1$.
The iteration ends when either $\xi_m=1$ or $\xi_m$ is a limit ordinal.
If $\xi_m > 1$ is a successor ordinal for every $m \in \dN$, the iteration never ends.
\end{enumerate}

Note that because $\xi$ is a limit ordinal we necessarily have $\xi_m < \xi$,
for every $m \in \dN$.

The next lemma states that odd iterations are finite.
\begin{lemma}
There is $m \in \dN$ such that either $\xi_m=1$ or $\xi_m$ is a limit ordinal.
\end{lemma}

\begin{proof}
Assume that the algorithm never terminates: $\xi_m > 1$ is a successor ordinal for every $m \in \dN$,
so that the algorithm generates an infinite sequence
$(h_1,\xi_1,a_1,h_2,\xi_2,a_2,\ldots)$.

We first argue that for every $m \in \dN$ one has
\begin{equation}
H_{\xi_m-1}(h_{m+1}) \supseteq H_{\xi_{m+2}-1}(h_{m+3}).
\end{equation}
Note that because $\xi_m$ and $\xi_{m+2}$ are successor ordinals,
the two ordinals $\xi_m-1$ and $\xi_{m+2}-1$ are well defined.
Now,
let $p$ be any play in $H_{\xi_{m+2}-1}(h_{m+3})$.
By Lemma \ref{lemma simple0}, $p \in H_{\xi_{m+2}}(h_{m+2})$.
Because $\xi_{m+1}-1 \leq \xi_{m+2}$, by the induction hypothesis (Part 4),
$p \in H_{\xi_{m+1}-1}(h_{m+2})$.
By Lemma \ref{lemma simple0}, $p \in H_{\xi_{m+1}}(h_{m+1})$.
Because $\xi_{m}-1 \leq \xi_{m+1}$, by the induction hypothesis (Part 4),
$p \in H_{\xi_{m}-1}(h_{m+1})$, as desired.

Because both $\xi_m$ and $\xi_{m+2}$ are successor ordinals,
Eq. (\ref{equ102}) implies that for every $m \in \dN$
\begin{equation}
\alpha_{\xi_m}(h_{m}) \leq \alpha_{\xi_{m+2}}(h_{m+2}).
\end{equation}
Because $\alpha_\xi(h_m) = \alpha_{\xi_m}(h_m)$ and
$\alpha_\xi(h_{m+2}) = \alpha_{\xi_{m+2}}(h_{m+2})$ we deduce that for every $m \in \dN$
\begin{equation}
\label{equ12}
\alpha_\xi(h_m) \leq \alpha_\xi(h_{m+2}).
\end{equation}

Because the payoffs are discrete, the
inequality $\alpha_\xi(h_m) \leq \alpha_\xi(h_{m+2})$ can be strict only finitely many times.
That is, there is $M \in \dN$ sufficiently large such that $\alpha_\xi(h_m) = \alpha_\xi(h_{m+2})$ for
every $m \geq M$.
Following the argument that Eq. (\ref{equ12}) holds we deduce that for each such $m$ we have
\begin{eqnarray}
\alpha_\xi(h_m) = \alpha_{\xi_m}(h_m) =
\alpha_{\xi_{m+1}-1}(h_{m+2}) = \alpha_{\xi_{m+2}}(h_{m+2}) = \alpha_\xi(h_{m+2}),
\end{eqnarray}
so that $\xi_{m+2} = \xi_{m+1}-1$.
Because this equality holds for every $m$ sufficiently large,
and because there is no infinite decreasing sequence of ordinals,
there is $m$ such that either $\xi_m=1$ or $\xi_m$ is a limit ordinal, as desired.
\end{proof}

\noindent\textbf{Even iterations:}

Let $h_1$ be the history that was generated by the previous iterations.
Then it is the output of the previous odd iteration.
Denote by $\lambda$ the last ordinal $\xi_m$ generated in the previous odd iteration.
Then in particular either $\lambda=1$ or $\lambda$ is a limit ordinal, and $\lambda < \xi$.
Moreover, $\alpha_\xi(h_1) = \alpha_\lambda(h_1)$.

By the induction hypothesis (Part 1 or Part 7),
there is a play $p \in H_\lambda(h_1)$ that is $\alpha_\lambda$-monotonic at $h_1$.

By the induction hypothesis (Parts 2, 3 and 4),
$\alpha_\xi(h') \geq \alpha_\lambda(h')$ for every prefix $h'$ of $p$ that extends $h_1$.
If $\alpha_\xi(h') = \alpha_\lambda(h')$ for every prefix $h'$ of $p$ that extends $h_1$,
the even iteration is infinite.
Otherwise,
the output of the current even iteration is the
shortest prefix $h'$ of $p$ that extends $h_1$
for which $\alpha_\xi(h') > \alpha_\lambda(h')$.

Because $p \in H_\lambda(h_1)$, the play $p$ is $\alpha_\lambda$-viable at $h_1$.
Because $\alpha_\xi(h') = \alpha_\lambda(h')$ for prefixes $h'$ of $p$ that extend $h_1$,
$p$ satisfies the condition of $\alpha_\xi$-viability for all such prefixes.

\begin{lemma}
\label{lemma6}
Let $p_*$ be the play that we just generated.
Then $p_*$ is $\alpha_\xi$-monotonic at $h$.
\end{lemma}

\begin{proof}
For the part of the play added in odd iterations the monotonicity was proved in (\ref{equ12}).
For the part of the play added in even iterations it follows from the construction.
\end{proof}

We are now ready to prove Part 7.
\begin{lemma}
$p_* \in H_\xi(h)$.
\end{lemma}

\begin{proof}
If the number of iterations if finite, so that the last even iteration is infinite,
then the play $p_*$ is $\alpha_\xi$-viable at $h$.
Indeed, by construction it is $\alpha_\xi$-viable at $h_0$, the history at the beginning of the last even iteration.
The claim now follows from Lemma \ref{lemma6}.

We now show that it cannot be that the number of iterations is infinite.
Denote by $(h_n)_{n \in \dN}$ all finite prefixes of $p_*$.
Because $p_*$ is $\alpha_\xi$-monotonic at $h$,
the two sequences $(\alpha(h_{2k}))_{h \preceq h_{2k}}$ and  $(\alpha(h_{2k+1}))_{h \preceq h_{2k+1}}$
are non-decreasing.
Because the range of the payoffs is finite, these two sequences are eventually constant.
At the beginning of an even iteration we have $\alpha_\xi(h_1) = \alpha_\lambda(h_1)$,
and the sequence $(\alpha_\lambda(h'))_{h'}$ is non-decreasing in the part of $p^*$ that is added in an even iteration
(for the definition of $\lambda$, see the construction for even iterations).
But $\alpha_\xi(h') \geq \alpha_\lambda(h')$ for every prefix $h'$ of $p_*$ that is added in the even iteration,
so that $\alpha_\xi(h') = \alpha_\lambda(h')$ for every such $h'$.
In particular, the even iteration does not end.
\end{proof}
\end{proof}

By Theorem \ref{theorem17}(2) and \ref{theorem17}(3) it follows that
\begin{theorem}
There is an ordinal $\xi_*$ such that
$\alpha_{\xi_*}(h)=\alpha_{\xi_*+1}(h)$ and $H_{\xi_*}(h)=H_{\xi_*+1}(h)$ for every $h \in H^1 \cup H^2$.
\end{theorem}

\subsection{Proof of Theorem \ref{theorem5}}
\label{subsection proof}

We now construct a pair $(\sigma^1_*,\sigma^2_*)$ of strategies, and show that they form
a subgame-perfect 0-equilibrium.

For every finite history $h$ choose an $\alpha_{\xi_*}$-viable play $p(h)$ that extends $h$ and that satisfies
\begin{equation}
u^{-i(h)}(p(h))  =\min_{p\in H_{\xi_*}(h)} u^{-i(h)}(p).
\end{equation}
If player $i$ deviates, and $h$ is the finite history right after the deviation
(so that $i = -i(h)$),
then $p(h)$ is an $\alpha_{\xi_*}$-viable play at $h$ that minimizes player's $i$'s payoff.

Let $\sigma^1_*$ be the following strategy:
Follow the play $p(\emptyset)$ as long as player 2 follows $p(\emptyset)$.
Suppose that at stage $2k_1$ player 2 deviates from $p(\emptyset)$.
From stage $2k_1+1$ and on follow the play $p(h_{2k_1+1})$ as long as player 2 follows this play.
Suppose that at stage $2k_2$ player 2 deviates from $p(h_{2k_1+1})$.
From stage $2k_2+1$ and on follow the play $p(h_{2k_2+1})$ as long as player 2 follows this play.
Continue this way.

The strategy $\sigma^2_*$ of player 2 is defined symmetrically.

We now show that $(\sigma^1_*,\sigma^2_*)$ is a subgame-perfect $0$-equilibrium.
To this end we fix a finite history $h \in H^1$
and we show that
\[ u^2(\sigma^1_*,\sigma^2 \mid h) \leq u^2(\sigma^1_*,\sigma^2_* \mid h), \ \ \ \forall \sigma^2 \in \Sigma^2. \]
One can use similar argument to show that the analog inequality for player 1 holds as well,
so that $(\sigma^1_*,\sigma^2_*)$ is indeed a subgame-perfect $0$-equilibrium.

Let $\sigma^2 \in\Sigma^2$ be any strategy of player 2.
Let $p_* = p(\sigma^1_*,\sigma^2_* \mid h)$ be the play induced by $(\sigma^1_*,\sigma^2_*)$
given $h$.
This is the play that is generated if player 2 does not deviate.
Let $p = p(\sigma^1_*,\sigma^2 \mid h)$ be the play when player 2 deviates to $\sigma^2$.

Denote by $2k_1, 2k_2,\ldots$ the stages where $\sigma^2$ and $\sigma^2_*$ differ along $p$;
in those stages player 1 observes the deviations of player 2. The sequence $(2k_j)_j$ may be finite or infinite.
Denote by $p_j = p(\sigma^1_*,\sigma^2)$ the play that player 1 start to follow at stage $2k_j$,
for each $j$.

We complete the proof by showing that
\begin{equation}
\label{equ 199}
u^2(p) \leq u^2(p_*).
\end{equation}
It is sufficient to show that
\begin{equation}
\label{equ 200}
u^2(p_j) \leq u^2(p_*), \ \ \ \forall j.
\end{equation}
If $\sigma^2$ and $\sigma^2_*$ differ only finitely many times along $p$,
Eq. (\ref{equ 199}) follows from Eq. (\ref{equ 200}).
If $\sigma^2$ and $\sigma^2_*$ differ infinitely many times along $p$,
then the sequence $(p_j)_{j \in \dN}$ converges to $p$, so that
Eq. (\ref{equ 199}) follows from Eq. (\ref{equ 200}) and the lower-semi-continuity of $u^2$.
This is the only place in the proof where the lower-semi-continuity of the payoff functions is used.

The proof of (\ref{equ 200}) is by induction on $j$.
Because $p_*$ is $\alpha_{\xi_*}$-viable at $h$,
$u^2(p_1) = \alpha_{\xi_*}(h_{2k_1}) \leq u^2(p_*)$.
For every $j \geq 1$, because $p_j$ is $\alpha_{\xi_*}$-viable at $h_{2k_j}$,
$u^2(p_{j+1}) = \alpha_{\xi_*}(h_{2k_j}) \leq u^2(p_{j})$,
which is at most $u^2(p_*)$ by the induction hypothesis.
The proof is now complete.

\subsection{A Folk Theorem}

Our construction enables us to characterize the set of plays that can arise in a subgame-perfect 0-equilibrium
in the game with discrete payoffs.

\begin{theorem}
\label{theorem2}
A play $p$ is induced by some subgame-perfect 0-equilibrium
if and only if $p$ is $\alpha_{\xi_*}$-viable.
\end{theorem}

\begin{proof}
If $p$ is $\alpha_{\xi_*}$-viable,
then the construction in Section \ref{subsection proof} shows that it is the play that is induced by some
subgame-perfect 0-equilibrium.

To see that the converse is true, we show that if $(\sigma^1_*,\sigma^2_*)$ is a
subgame-perfect 0-equilibrium, then $p(\sigma^1_*,\sigma^2_* \mid h)$ is $\alpha_\xi$-viable at $h$,
for every ordinal $\xi$ and every finite history $h$.

Because $\alpha_1(h)$ is the value of the subgame that starts at $h$,
and because every subgame-perfect 0-equilibrium is an equilibrium in the subgame that starts at $h$,
the claim follows for  $\xi=1$.

Suppose now that the claim holds for an ordinal $\xi$.
Let $h$ be any finite history.
Because the claim holds for $\xi$,
the play $p(\sigma^1_*,\sigma^2_* \mid h)$ is $\alpha_\xi$-viable at $(h,a)$ for every $a$,
and therefore in any subgame-perfect 0-equilibrium of the subgame that starts at $h$,
the payoff to player $i(h)$ is at least $\alpha_{\xi+1}(h)$.
This implies that $p(\sigma^1_*,\sigma^2_* \mid h)$ is $\alpha_{\xi+1}$-viable at $h$,
for every $h$.

Finally, the definition of $\alpha_\xi$ for limit ordinals $\xi$ implies
that if $p(\sigma^1_*,\sigma^2_* \mid h)$ is $\alpha_\lambda$-viable at $h$ for every ordinal $\lambda < \xi$,
then it is also $\alpha_\xi$-viable at $h$.
\end{proof}

\section{Comments}
\label{comments}

\subsection{Chance moves}

Borel games are deterministic, and the sequence of actions chosen by the players
uniquely determines the outcome.
In many situations there are chance moves along the game,
where actions are chosen according to a known probability distribution.
This situation is equivalent to the case where there is a third player who follows a specific
non-deterministic strategy,
whatever the other players play.
Our proof can be adapted to this more general situation,
and this will be done elsewhere.

\subsection{Positive recursive Borel games}

Recursive Borel games are games where some finite histories are terminating,
in the sense that once they occur the payoff is determined
(and the play that follows them does not affect the players' payoffs),
and the payoff of every infinite (non-terminating) play is 0.
Many positional games that are studied in the computer science literature have this form.
The significance of this class of games to game theory
was exhibited in the context of stochastic games by Vieille (2000),
who used it as a step to proving the existence of an equilibrium payoff in every two-player stochastic game.
A recursive Borel game is called positive if both $u^1$ and $u^2$ are positive functions.

Flesch et al. (2008) studied positive recursive Borel game with finitely many states;
these are positional games that are played on a finite directed graph, where each vertex is controlled by some player,
and when the game reaches some vertex, the controlling player can choose whether to terminate the game,
or whether to continue the game by choosing one of the edges that leaves the vertex.
The terminal payoff depends only on the vertex where termination occurred,
and not on the whole past play.

Flesch et al. (2008) prove that every such game admits
a subgame-perfect 0-equilibrium.%
\footnote{When transitions are random, Flesch et al. (2008) prove the existence of a
subgame-perfect $\ep$-equilibrium, for every $\ep > 0$.}
In their proof, they define for every vertex $s$ a sequence $(\alpha_k(s))_{k \in \dN}$
that is similar to our sequence $(\alpha_\xi(h))_\xi$,
they prove that this sequence is non-decreasing,
and, because there are finitely many vertices, they deduce that there is $k_*\in \dN$
such that $\alpha_{k_*+1}(s) = \alpha_{k_*}(s)$
for every vertex $s$.
They then use a similar construction of the subgame-perfect 0-equilibrium as the one that we used.

In Borel games every history is a different vertex.
Therefore one needs to employ a much more delicate construction,
that differs from the one in Flesch et al. (2008) in two respects.
First, when the number of vertices is infinite,
there need not by $k_*\in \dN$ such that $\alpha_{k_*+1}(s) = \alpha_{k_*}(s)$
for every vertex $s$,
and therefore $(\alpha_\xi(h))_\xi$ should be defined for every ordinal.
Second, one also needs to take into account the possible plays,
and introduce the sets $(H_\xi(h))_\xi$.
In fact, this paper grew from an attempt to understand why the proof in Flesch et al. (2008) fails for positive
recursive games with infinitely many vertices.

It turns out that for positive recursive Borel games the construction can be simplified,
and a single odd iteration is sufficient to show that $H_\xi(h)$ is not empty for limit ordinals $\xi$.

\subsection{Borel games with general payoffs}

Example 3 in Solan and Vieille (2003) shows that without the condition that payoffs are lower-semi-continuous,
a subgame-perfect $\ep$-equilibrium need not exist.
However, Solan and Vieille (2003)
show that a subgame-perfect $\ep$-equilibrium does exist if one allows behavior strategies.
The existence of a subgame-perfect $\ep$-equilibrium in behavior strategies
was proved in other setups where the payoff functions are not lower-semi-continuous,
see, e.g., Solan (2005) and Mashiah-Yaakovi (2008).

In our proof the lower-semi-continuity of the payoff functions was used only in the last part,
to show that any deviation $\sigma^2$ that differs from $\sigma^2_*$ infinitely many times cannot be profitable,
as soon as any deviation $\sigma^2$ that differs from $\sigma^2_*$ finitely many times is not profitable.
We do not know how the proof should be adapted to handle general payoff functions.

In fact, the following example shows that our definition of $\alpha_\xi$ and $H_\xi$ is not appropriate for
general Borel games.
Consider a Borel game with $A = \{a,b\}$.
The payoff functions of the two players are as follows:
\[
\begin{array}{l|l|l}
\hbox{Condition} & u^1(h) & u^2(h) \\
\hline
\hbox{Both players played } b \hbox{ finitely many times} & 2 & 2\\
\hbox{Only player 1 played } b \hbox{ finitely many times} & 2 & 1\\
\hbox{Only player 2 played } b \hbox{ finitely many times} & 1 & 2\\
\hbox{No player played } b \hbox{ finitely many times} & 0 & 0
\end{array}
\]
Note that $u^1$ and $u^2$ are not lower-semi-continuous.
Playing $b$ finitely many times is a dominant strategy for both players,
so that the unique subgame-perfect 0-equilibrium payoff is $(2,2)$.
However, one can verify that $\alpha_\xi(h) = 1$ for every finite history $h$ and every ordinal $\xi$,
so that the folk theorem, Theorem \ref{theorem2}, does not hold,
and our construction of the subgame-perfect 0-equilibrium in the proof of Theorem \ref{theorem01} is invalid.

\end{document}